\documentclass[12pt]{amsart}
\usepackage{amssymb, amsmath}
\newtheorem{question}{Question}[section]
\newtheorem{theorem}[question]{Theorem}
\newtheorem{lemma}[question]{Lemma}
\newtheorem{corollary}[question]{Corollary}
\newtheorem{example}[question]{Example}
\newtheorem{fact}[question]{Fact}

\newtheorem{definition}[question]{Definition}

\title{A note on discrete sets}
\author{Santi Spadaro}

\address{
Department of Mathematics and Statistics\\
Auburn University\\
221 Parker Hall\\
Auburn, Alabama -- USA 36849-5310}

\email{spadasa@auburn.edu}

\thanks{Research partially supported by National Science Foundation grant DMS-0405216 (Principal Investigator - Dr. Gary Gruenhage)}

\subjclass[2000]{54A25}

\keywords{discrete set, dispersion character, compact space, Eberlein compact, free sequence, elementary submodel}

\begin{document}
\baselineskip.525cm

\begin{abstract}
We give several partial positive answers to a question of Juh\'asz and Szentmikl\'ossy regarding the minimum number of discrete sets required to cover a compact space. We study the relationship between the size of discrete sets, free sequences and their closures with the cardinality of a Hausdorff space, improving known results in the literature.
\end{abstract}

\maketitle

\section{Introduction}
How many discrete sets does it take to cover a compact space? Do discrete sets reflect the cardinality of a compact space? These questions have been considered by many authors.

Call $dis(X)$ the least number of discrete sets required to cover $X$. Gruenhage (\cite{G3}) proved that $dis(X) \geq \mathfrak{c}$ for every compact space without isolated points, thus answering a question of Juh\'asz and van Mill \cite{JM}. It is unknown whether compactness can be replaced by countable compactness.
By exploiting a lemma of Gruenhage, and yet using a completely different approach, Juh\'asz and Szentmikl\'ossy \cite{JS} proved that in every compact space with $\chi(x,X) \geq \kappa$ for every $x \in X$ we have $dis(X) \geq 2^\kappa$, thus generalizing both Gruenhage's theorem and the classical \v Cech-Pospi\v sil theorem. 

Call $\Delta(X)$ the least cardinality of a non-empty open set in $X$. Since in every compact space where every point has character at least $\kappa$ we have $\Delta(X) \geq 2^\kappa$, Juh\'asz and Szentmikl\'ossy naturally ask the following question.

\begin{question} \label{q1}
\cite{JS} Is $dis(X) \geq \Delta(X)$ for every compact space $X$?
\end{question}

In the first part of this note we will give several partial positive answers to the previous question. 

We consider special classes of compact spaces (compact $T_5$ spaces, compact LOTS, polyadic compacta, Gul'ko compacta...) and prove that Juh\'asz and Szentmikl\'ossy's inequality is true for them. A few results outside of the compact realm are obtained as a byproduct, for example we determine the least number of discrete sets required to cover a $\Sigma$-product. Also, sometimes we can replace compactness by a weaker property (for example, the Baire property).

Let now $g(X)=\sup \{|\overline{D}|: D \subset X$ is discrete. Alas, Tkachuk and Wilson \cite{ATW} ask whether $g(X) \leq \mathfrak{c}$ implies that $|X| \leq \mathfrak{c}$ for every compact $X$. Only consistent negative answers are known for this question (see \cite{D2}). On the other hand, Alas provides a partial consistent positive answer in the following theorem.

\begin{theorem}
\cite{A} (MA) If $X$ is compact, $\hat{s}(X) \leq \mathfrak{c}$ and $g(X) \leq \mathfrak{c}$ then $|X| \leq \mathfrak{c}$.
\end{theorem}

The condition $\hat{s}(X) \leq \mathfrak{c}$ just means that every discrete set in $X$ has size $<\mathfrak{c}$. Another partial positive answer is provided by the following theorem of Alan Dow.

\begin{theorem}
\cite{D2} If $X$ is a compact space of countable tightness such that $g(X) \leq \mathfrak{c}$ then $|X| \leq \mathfrak{c}$.
\end{theorem}

In the second part of this note we are going to prove a common generalization of the above results that takes them out of the compact realm. Further investigations on when cardinality is reflected by discrete sets or even free sequences will follow.

All spaces are assumed to be Hausdorff. A space is called \emph{crowded} if it has no isolated points. All undefined notions can be found in \cite{E} and \cite{J}. The spread, cellularity, tightness, weight, $\pi$-character and the number of regular open sets of $X$ will be denoted respectively by $s(X)$, $c(X)$, $t(X)$, $w(X)$, $\pi \chi(X)$ and $\rho(X)$.

\section{Covering a compact space by discrete sets}

Testing a conjecture about compact spaces on compact hereditarily normal spaces is quite a natural thing to try, and indeed, Juh\'asz and Van Mill already did that for the inequality $dis(X) \geq \mathfrak{c}$, before Gruenhage proved it to be true for every compact Hausdorff space. 

\begin{theorem} \label{perfect} (\cite{G3})
Let $f: X \to Y$ be a perfect map. Then $dis(X) \geq dis(Y)$.
\end{theorem}

Let $\kappa^\omega$ be the product of countably many copies of the discrete space $\kappa$. The following was proved in \cite{S}.

\begin{lemma} \label{metr}
For every Baire metric space $dis(X) \geq \Delta(X)$. In particular $dis(\kappa^\omega) = \kappa^\omega$.
\end{lemma}

A \emph{cellular} family is a family of pairwise disjoint open sets in $X$. The following lemma is crucial to most of our results in this section.

\begin{lemma} \label{cell}
Let $X$ be a compact space whose every open set contains a cellular family of cardinality $\kappa$. Then $dis(X) \geq \kappa^\omega$.
\end{lemma}

\begin{proof}
Use regularity of $X$ to find a cellular family $\{U_\alpha: \alpha < \kappa \}$ such that the closures of its members are pairwise disjoint. Suppose you have constructed open sets $\{U_\sigma : \sigma \in \kappa^{<n} \}$. Then let $\{U_{\sigma^{\frown} \alpha} : \alpha \in \kappa \}$ be a cellular family inside $U_\sigma$ such that the closures of its members are pairwise disjoint and contained in $U_\sigma$.

For each $f \in \kappa^\omega$ let $F_f=\bigcap_{n \in \omega} \overline{U_{f \upharpoonright n}}$, which is a non-empty set because of compactness, and set $Z=\bigcup_{f \in \kappa^\omega} F_f$.

We are now going to show a perfect map $\Phi$ from $Z$ onto $\kappa^\omega$. By Theorem $\ref{perfect}$ and Lemma $\ref{metr}$ we will get that $dis(X) \geq \kappa^\omega$.

Define $\Phi$ simply as $\Phi(x)=f$ whenever $x \in F_f$. It is easy to see that the $F_f$'s are pairwise disjoint, so $\Phi$ is well-defined. Moreover, $\Phi$ is clearly continuous, onto and has compact fibers.

The following characterization of closed maps is well-known (see \cite{E}, Theorem 1.4.13).

\begin{fact}
A mapping $f: X \to Y$ is closed if and only if for every point $y \in Y$ and every open set $U \subset X$ which contains $f^{-1} (y)$, there exists in $Y$ a neighbourhood $V$ of the point $y$ such that $f^{-1}(V) \subset U$.
\end{fact}

Let now $f \in \kappa^\omega$, and $U$ be an open set in $Z$ such that $\Phi^{-1}(f)=F_f=\bigcap_{n \in \omega} \overline{U_{f \upharpoonright n}} \subset U$. By compactness, we can find an increasing sequence of integers $\{j_k: 1 \leq k \leq n \}$ such that $\overline{U_{f \upharpoonright j_n}}=\bigcap_{1 \leq k \leq n} \overline{U_{f \upharpoonright j_k}} \subset U$.

So let $B(f \upharpoonright j_n)$ be the basic neighbourhood in $\kappa^\omega$ determined by $f \upharpoonright j_n$. Then $\Phi^{-1} (B(f \upharpoonright j_n)) \subset U_{f \upharpoonright j_n} \subset U$, which proves $\Phi$ is closed. 
\end{proof}

\begin{theorem} \label{coll}
Let $X$ be an hereditarily collectionwise Hausdorff compact space. Then $dis(X) \geq \Delta(X)$.
\end{theorem}

\begin{proof}
Recall that cellularity and spread coincide for hereditarily collectionwise Hausdorff spaces (see \cite{J}, 2.23 a)). So if $c(G) < \Delta(X)$, for some open set $G \subset X$ we also have $s(G) < \Delta(X) \leq \Delta(G)$. Hence $dis(X) \geq \Delta(X)$. 

Suppose now that $c(G) \geq \Delta(X)$ for every open set $G \subset X$. If $\Delta(X)$ is a successor cardinal then every open set contains a cellular family of size $\Delta(X)$, and hence, in view of Lemma $\ref{cell}$ we have $dis(X) \geq \Delta(X)$. 

If $\Delta(X)$ is a limit cardinal then, again by Lemma $\ref{cell}$, every open set contains a cellular family of size $\kappa$ for every $\kappa < \Delta(X)$. Hence $dis(X) \geq \kappa$ for every $\kappa < \Delta(X)$, which implies $dis(X) \geq \Delta(X)$ again.
\end{proof}

\begin{corollary}
For every compact LOTS $X$, $dis(X) \geq \Delta(X)$.
\end{corollary}

\begin{proof}
Compact LOTS are monotonically normal, which implies collectionwise normal, and monotone normality is hereditary (see \cite{G1}).
\end{proof}

From Theorem $\ref{coll}$ it also follows that, under V=L, $dis(X) \geq \Delta(X)$ for every compact hereditarily normal space $X$. Indeed, Stephen Watson \cite{W} proved that compact $T_5$ spaces are hereditarily collectionwise Hausdorff in the constructible universe. We can do better, and prove that $dis(X) \geq \Delta(X)$ for $X$ compact $T_5$ under a slight weakening of GCH.

\begin{theorem}
(for every cardinal $\kappa$, $2^\kappa < 2^{\kappa^+}$) Let X be a compact $T_5$ space. Then $dis(X) \geq \Delta(X)$.
\end{theorem}

\begin{proof}
Suppose first that $c(G) < \Delta(X)$ for some open set $G$. Since $c(G)=c(\overline{G})$ and $\overline{G}$ is compact $T_5$ we can assume that $X=\overline{G}$.

Let $\kappa=c(X)$. By Shapirovskii's bound on the number of regular open sets (see \cite{J}, 3.21) we have $\rho(X) \leq 2^\kappa$. Note that $\kappa^+ \leq \Delta(X)$. If $dis(X) < \Delta(X)$ then we would have $s(X) \geq \Delta(X)$ and hence we could find a discrete $D \subset X$ such that $|D| \geq \kappa^+$. By Jones' Lemma (see \cite{J}, 3.1), $\rho(X) \geq 2^{\kappa^+}> 2^\kappa$, which contradicts our upper bound for the number of regular open sets.

If $c(G) \geq \Delta(X)$ for every open set $G$, then reasoning as in the last few lines of the proof of Theorem $\ref{coll}$ we can conclude that $dis(X) \geq \Delta(X)$.
\end{proof}

\begin{question}
Is it true in ZFC that $dis(X) \geq \Delta(X)$ for every compact $T_5$ space?
\end{question}

A trivial observation is that all compact metrizable spaces satisfy $dis(X) \geq \Delta(X)$.

The two most popular generalizations of compact metrizable spaces are dyadic compacta and Eberlein compacta. In fact, they are two somewhat opposite classes, as their intersection is precisely the class of compact metrizable spaces (see \cite{Ar}). 

This made us wonder whether $dis(X) \geq \Delta(X)$ was true for them. In fact, we are able to prove that for the weaker classes of \emph{polyadic} and \emph{Gul'ko} compacta. To achieve that we first need to prove that $dis(X)$ is always bounded below by the tightness. Recall that a space is called initially $\kappa$-compact if every set of cardinality $ \leq \kappa$ has a complete accumulation point.

\begin{lemma} (\cite{GJS}) \label{lemgjs}
Let $X$ be an initially $\kappa$-compact space such that $dis(X) \leq \kappa$. Then $X$ is compact.
\end{lemma}

\begin{lemma} \label{tight}
If $X$ is compact then $dis(X) \geq t(X)$.
\end{lemma}

\begin{proof}
Suppose by contradiction that $\kappa=dis(X) < t(X)$. Let $A \subset X$ be a non-closed set, and $[A]_\kappa$ be its $\kappa$-closure, that is, the union of the closures of its subsets of cardinality $\kappa$. If we could prove that this last set is closed then we would have $t(X) \leq \kappa$, which is what we want.

If $[A]_\kappa$ is not closed then it cannot be initially $\kappa$-compact, or otherwise, since $dis([A]_\kappa) \leq \kappa$, it would be compact by Lemma $\ref{lemgjs}$. So there is $B \subset [A]_\kappa$ such that $|B| \leq \kappa$ and $B$ has no point of complete accumulation in $[A]_\kappa$; then, by compactness, there is a point $x \notin [A]_\kappa$ that is of complete accumulation for $B$. But this contradicts the well-known and easy to prove fact that $[[A]_\kappa]_\kappa=[A]_\kappa$.

\end{proof}

A compactum is called \emph{polyadic} if it is the continuous image of some power of the one-point compactification of some discrete set.

The following lemmas are due to Gerlits.

\begin{lemma} 
\cite{Ge1} Let $X$ be polyadic and $A \subset X$. Then there is a polyadic $P \subset X$ such that $A \subset P$ and $c(P) \leq c(A)$.
\end{lemma}

\begin{lemma}
\cite{Ge2} If $X$ is polyadic then $w(X)=t(X) \cdot c(X)$.
\end{lemma}

\begin{theorem}
For a polyadic compactum $X$ we have $dis(X) \geq \Delta(X)$.
\end{theorem}

\begin{proof}
If $c(U) \geq \Delta(X)$ for any open set $U \subset X$ then we are done by Lemma $\ref{cell}$. If there exists some open $U$ such that $c(U) < \Delta(X)$, then let $P$ be a polyadic space such that $U \subset P$ and $c(P) \leq c(U)$. Assume $dis(P) < \Delta(X)$. Then $t(P) < \Delta(X)$, which implies $s(P) \leq w(P) < \Delta(X)$, and we are done, since $|P| \geq \Delta(X)$.
\end{proof}

Recall that an \emph{Eberlein compactum} is a compact space which embeds in $C_p(Y)$ for some compact $Y$. Equivalently, a space is an Eberlein compactum if and only if it is a weakly compact subspace of a Banach space. A \emph{Gul'ko compactum} is a compact space $X$ such that $C_p(X)$ is a Lindel\"of $\Sigma$-space. A \emph{Corson compactum} is a compact space with embeds in a $\Sigma$-product of lines. The following chain of implications holds.

$$\mbox{Eberlein} \Rightarrow \mbox{Gul'ko} \Rightarrow \mbox{Corson}$$

Recall that a space is called \emph{meta-Lindel\"of} if every open cover has a point-countable open refinement. 

\begin{lemma}
Let $X$ be an hereditarily meta-Lindel\"of space such that $dis(X) \leq \kappa$. If $A \subset X$ is such that $|A| \leq \kappa$ then $|\overline{A}| \leq \kappa$.
\end{lemma}

\begin{proof}

If $\kappa < \omega$ then the statement is obviously true. Assume that $\kappa$ is infinite, and let $X=\bigcup_{\alpha < \kappa} D_\alpha$, where each $D_\alpha$ is discrete. Let $B_\alpha = \overline{A} \cap D_\alpha$. For every $x \in B_\alpha$, let $U_x$ be an open set such that $U_x \cap B_\alpha=\{x\}$. Then $\bigcup_{x \in B_\alpha} U_x$ is meta-Lindel\"of, and hence $\{U_x : x \in B_\alpha \}$ has a point-countable open refinement $\mathcal{V}_\alpha$. Now for every $x \in B_\alpha$ choose $V_x \in \mathcal{V}_\alpha$ such that $x \in V_x$ and let $\mathcal{U}_\alpha=\{V_x : x \in B_\alpha \}$. Clearly $|\mathcal{U}_\alpha|=|B_\alpha|$ and for all $U \in \mathcal{U}_\alpha$, $U \cap A \neq \emptyset$. Fix some well-ordering of $A$ and define a function $f : \mathcal{U_\alpha} \to A$ by:

$$f(U)=\min \{a \in A : a \in U \}.$$ 

Point-countability of $\mathcal{U}_\alpha$ implies that $|f^{-1}(a)| \leq \aleph_0$ for every $a \in A$, and therefore $|B_\alpha|=|\mathcal{U}_\alpha| \leq |A| \cdot \aleph_0 \leq \kappa$.

Since $\overline{A}=\bigcup_{\alpha \in \kappa} B_\alpha$ it follows that $|\overline{A}| \leq \kappa$.
\end{proof}

\begin{theorem}
Let $X$ be an hereditarily meta-Lindel\"of space containing a dense Baire metrizable subset. Then $dis(X) \geq \Delta(X)$.
\end{theorem}

\begin{proof}
Let $M \subset X$ be a dense metrizable subset and suppose by contradiction that $dis(X) < \Delta(X)$. Then, by the previous lemma we have $\Delta(M)=\Delta(X)$. So $dis(X) \geq dis(M) \geq \Delta(M)=\Delta(X)$, which is a contradiction.
\end{proof}

\begin{corollary}
For every Gul'ko compactum $X$ we have $dis(X) \geq \Delta(X)$.
\end{corollary}

\begin{proof}
Yakovlev (\cite{Y}) proved that every Corson compactum is hereditarily meta-Lindel\"of and Gruenhage (\cite{G2}) proved that every Gul'ko compactum contains a dense Baire metrizable subset.
\end{proof}

We are sorry to admit that we haven't been able to answer the following two questions.

\begin{question}
Is it true in ZFC that $dis(X) \geq \Delta(X)$ for every Corson compact $X$?
\end{question}

\begin{question}
Is $dis(X) \geq \Delta(X)$ for every compact space with a (Baire) dense metrizable subset?
\end{question}

As an application of the results in this section we are now going to determine how many discrete sets are needed to cover the $\Sigma$-product of a Cantor cube.

\begin{theorem} \label{sigma}
$dis(\Sigma(2^{\kappa})) = \kappa^\omega$.
\end{theorem}

To prove that we will embed in $\Sigma(2^{\kappa})$ an Eberlein compactum $X$ for which $\Delta(X)=\kappa^\omega$. 

Recall that a family $\mathcal{A}$ of subsets of a set $T$ is called \emph{adequate} if:

\begin{enumerate}
\item For every $A \in \mathcal{A}$, $\mathcal{P}(A) \subset \mathcal{A}$.
\item If $[A]^{<\omega} \subset \mathcal{A}$ then $A \in \mathcal{A}$.
\end{enumerate}

It is easy to see that $\mathcal{A}$ with the topology inherited from the product space $2^T$ is closed, and hence compact. Such a space is called an \emph{adequate compactum}. Adequate families are one of the most useful tools for constructing Corson compacta: especially handy is the adequate family of all chains of a partial order. If the partial order has no uncountable chains, then the corresponding adequate compactum is Corson. 

Leiderman and Sokolov characterized all adequate Eberlein compacta.

\begin{theorem} \label{leisok} (\cite{LS})
Let $X$ be an adequate compact embedded in $2^T$. Then $X$ is an Eberlein compact if and only if there is a partition $T=\bigcup_{i \in \omega} T_i$ such that $|supp(x) \cap T_i|<\aleph_0$ for each $x \in X$ and $i \in \omega$.
\end{theorem} 

The next example is a modification of an example due to Leiderman and Sokolov. Their original space was a strong Eberlein compactum (a weakly compact subset of a Hilbert space), and hence scattered. Our space is far from being scattered.

\begin{example}
Let $\kappa$ be any infinite cardinal. There is an Eberlein compactum, embedded in $2^\kappa$, such that $\Delta(X)=\kappa^\omega$.
\end{example}

\begin{proof}
Let $W_0=Lim(\kappa)$ and let $\{x_\alpha : \alpha \in \kappa \}$ be an increasing enumeration of $W_0$. Let $W_i=\{x_\alpha+i : \alpha \in \kappa \}$. Now let $T=\bigcup_{i \in \omega} W_i \times (W_i \cup \{-i\})$. Define an order on $T$ as follows : $(\alpha_1, \beta_1) < (\alpha_2, \beta_2)$ if and only if $\alpha_1 < \alpha_2$ and $\beta_1 > \beta_2$. Then every chain in $T$ is countable, so the adequate compact $X$ constructed from the adequate family consisting of all chains in $T$ is Corson. Moreover, the partition in the definition of $T$, along with Theorem $\ref{leisok}$ shows that $X$ is Eberlein. It remains to check that $\Delta(X)=\kappa^\omega$. To see that, let $U$ be any basic open set. Then $U$ is the set of all chains containing some fixed finite chain $\{(\alpha_i, \beta_i) : i \leq k \}$, enumerated in increasing order, and not containing a finite number of fixed elements. Let $j \in \omega$ be such that $\alpha_k \in W_j$. Now, for all but finitely many increasing chains $\{\alpha_{j+k} : k \geq 1 \}$, with $\alpha_{j+k} \in W_{j+k}$ for every $k \geq 1$ and $\alpha_j+1 > \alpha_j$ we have that $\{(\alpha_i, \beta_i) : i \leq k \} \cup \{(\alpha_{j+k}, -(j+k)) : k \geq 1 \} \in U$. Now the set of all such chains has cardinality $\kappa^\omega$, since there is a natural bijection between that set and the set of all countable increasing sequences in $\kappa$. 
\end{proof}

Every $\Sigma$-product of compact spaces is countably compact, which reminds us of the following question.

\begin{question}
Is $dis(X) \geq \mathfrak{c}$ for $X$ countably compact crowded?
\end{question}

The starting point for our next pair of results is the following easy observation.

\begin{theorem}
Let $X$ be a homogeneous compactum. Then $dis(X) \geq \Delta(X)$.
\end{theorem}

\begin{proof}
Combining Arhangel'skii's theorem with the Juh\'asz-Szentmikl\'ossy's result cited in the introduction we get
$dis(X) \geq 2^{\chi(X)} \geq \Delta(X)$.
\end{proof}

A space is \emph{homogeneous with respect to character} if $\chi(x,X)=\chi(y,X)$ for any $x,y \in X$. A space $X$ is \emph{power homogeneous} if $X^\kappa$ is homogeneous for some $\kappa$. 

The following lemma is due to Juh\'asz and Van Mill.

\begin{lemma} \label{char} (\cite{JM})
Every infinite compactum contains a point $x$ with $\chi(x,X) < dis(X)$.
\end{lemma}

We are also going to need a couple of results from Guit Jan Ridderbos' PhD Thesis. 

\begin{lemma} \label{rid} (\cite{R})
Let $X$ be power homogeneous. If the set of all points of $\pi$-character $\kappa$ is dense in $X$, then $\pi \chi(X) \leq \kappa$.
\end{lemma}

\begin{lemma} \label{rid2} (\cite{R})
Let $X$ be a power-homogeneous space of pointwise countable type such that $\pi \chi(X) \leq \kappa$. Then either $\chi(X) \leq \kappa$ or $X$ is homogeneous with respect to character.
\end{lemma}

We are going to prove that under the GCH if a power-homogeneous compactum is \emph{not too big} then it satisfies Juh\'asz and Szentmikl\'ossy's inequality. We need the following lemma, which, in a sense, says that the gap between $\Delta(X)$ and $dis(X)$ can't be too large for power-homogeneous compacta.

\begin{lemma} \label{log}
Let $X$ be a power homogeneous compactum. Then $\Delta(X) \leq 2^{dis(X)}$.
\end{lemma}

\begin{proof}
Suppose by way of contradiction that $dis(X) \leq \kappa$ but $|U| > 2^\kappa$ for every open $U \subset X$. Then by Lemma $\ref{char}$ the set of all points of character less than $\kappa$ is dense $X$, which implies $\pi \chi(X) \leq \kappa$. If $\chi(X) \leq \kappa$, then, by Arhangel'skii's Theorem, $|X| \leq 2^\kappa$, which contradicts our initial assumption. 
Otherwise $\chi(X) \geq \kappa^+$ and $X$ is homogeneous with respect to character, which even implies $dis(X) \geq 2^{\kappa^+}$, again a contradiction.
\end{proof}

\begin{theorem}
(GCH) Let $X$ be a power-homogeneous compactum. Then $dis(X) \geq \min \{\Delta(X), \aleph_\omega\}$.
\end{theorem}

\begin{proof}
Suppose by contradiction that $dis(X) < \Delta(X)$ and $dis(X) < \aleph_\omega$. Then $dis(X)=\aleph_{n+1}$ for some $n \in \omega$. By Lemma $\ref{char}$ the space $X$ contains a dense set of $G_{\aleph_n}$ points, and hence $\pi \chi (X) \leq \aleph_n$ by Lemma $\ref{rid}$. If $X$ were homogeneous with respect to character then $dis(X) \geq \Delta(X)$ and we would get a contradiction. So, by Lemma $\ref{rid2}$, $\chi(X) \leq \aleph_n$ and hence $\Delta(X) \leq 2^{\aleph_n}=\aleph_{n+1}$. Now by Lemma $\ref{log}$ and GCH we have that $dis(X)^+=\Delta(X)$. So $\Delta(X)=\aleph_{n+2}$ and we get the desired contradiction.
\end{proof}

\begin{corollary}
(GCH) If $X$ is a power-homogeneous compactum such that $|X| \leq \aleph_\omega$ then $dis(X) \geq \Delta(X)$.
\end{corollary}

If $|X| \leq \aleph_3$ then we need only assume CH by a different proof.

\begin{theorem}
(CH) Let $X$ be a power-homogeneous compactum. Then $dis(X) \geq \min \{\Delta(X), \omega_3\}$.
\end{theorem}

\begin{proof}
Suppose that $\beta \omega$ does not embed in $X$, then $X$ does not map onto $I^{\omega_1}$ (see the proof of \cite{J}, 3.22) and hence, as a consequence of Shapirovskii's Theorem on maps onto Tychonoff cubes, the set of all points of countable $\pi$-character is dense in $X$. Therefore, by Lemma $\ref{rid}$, $\pi \chi(X) \leq \omega$. If $\chi(X) \leq \omega$, then $|X| \leq \omega_1$, by Arhangel'skii's theorem, and since $dis(X) \geq \omega_1$ holds for every compactum, we are done. Otherwise, $X$ is homogeneous with respect to character, and hence $|X| \leq 2^{\chi(X)} \leq dis(X)$, by Juh\'asz and Szentmikl\'ossy's result.

If $\beta \omega$ embeds in $X$ then $dis(X) \geq 2^{\omega_1}$. Suppose that $dis(X) < \omega_3$, that is $dis(X) \leq \omega_2$. Then, by Lemma $\ref{char}$, $X$ contains a dense set of $G_{\omega_1}$ points. If $\chi(X) \leq \omega_1$, then $\Delta(X) \leq 2^{\omega_1}$ and we are done. Otherwise, $X$ is homogeneous with respect to character, and $dis(X) \geq \Delta(X)$ is true again.
\end{proof}

\begin{corollary}
(CH) If $X$ is a power-homogeneous compactum such that $|X| \leq \omega_3$ then $dis(X) \geq \Delta(X)$.
\end{corollary}

\begin{question}
Is $dis(X) \geq \Delta(X)$ true for every power-homogeneous compactum?
\end{question}

\section{Closures of discrete sets and cardinality}

Alas, Tkachuk and Wilson \cite{ATW} asked whether a compact space in which the closure of every discrete set has size $\leq \mathfrak{c}$ must have size $\leq \mathfrak{c}$. 

In \cite{A} Ofelia Alas proves the following theorem, by way of a partial positive answer.

\begin{theorem}
(MA) Let $X$ be a Lindel\"of regular weakly discretely generated space such that $\hat{s}(X) \leq \mathfrak{c}$ and $|\overline{D}| \leq \mathfrak{c}$ for every discrete $D \subset X$. Then $|X| \leq \mathfrak{c}$.
\end{theorem}

We are going to prove that regular, Lindel\"of and weakly discretely generated can all be dropped from the above theorem. But, first of all let's define four cardinal functions that will be handy in our study of this and related problems. Recall that a sequence $\{x_\alpha: \alpha < \kappa \}$ is said to be \emph{free} if, for every $\gamma < \kappa$ we have $\overline{\{x_\alpha : \alpha \leq \gamma \}} \cap \overline{\{x_\alpha : \alpha > \gamma \}} = \emptyset$. Every free sequence is a discrete set.

\begin{definition}
Set $\hat{s}(X)=\min\{\kappa:$ if $A \subset X$ and $|A|=\kappa$ then $A$ is not a discrete set $\}$ and $\hat{F}(X)=\min\{\kappa:$ if $A \subset X$ and $|A|=\kappa$ then $A$ is not a free sequence $\}$.
\end{definition}

\begin{definition}
Set $g(X)=\sup \{|\overline{D}| : D \subset X$ is discrete $\}$ (the \emph{depth} of $X$) and $b(X)=\sup \{|\overline{F}|: F \subset X$ is a free sequence $\}$ (the \emph{breadth} of $X$).
\end{definition}

The condition $g(X) \leq \kappa$ appears to be a lot stronger than $b(X) \leq \kappa$. In fact, while the former implies that $|X| \leq 2^\kappa$ (simply observe that the hereditarily Lindel\"of number is discretely reflexive \cite{ATW} and use De Groot's inequality $|X| \leq 2^{hL(X)}$), the latter alone does not put any bound on the cardinality of $X$. For example, the one-point compactification of a discrete set of arbitrary cardinality satisfies $b(X)= \omega$.

Before proving our first theorem, we need a little lemma about elementary submodels, and an old lemma of Shapirovskii. All one needs to know about elementary submodels to read this section can be found in \cite{D1}. The following lemma is probably well-known. However, we include a proof of it anyway since we could not find a direct reference to it.

\begin{lemma}
Suppose $\mathfrak{c}$ is a regular cardinal. Let $\theta \geq (2^{<\mathfrak{c}})^+$ be a regular cardinal and $A \subset H(\theta)$ be a set of size $\leq 2^{<\mathfrak{c}}$. Then there is an elementary submodel $M \prec H(\theta)$ such that $A \subset M$, $|M| = 2^{<\mathfrak{c}}$ and $M$ is $\lambda$-closed for every $\lambda < \mathfrak{c}$.
\end{lemma}

\begin{proof}
It follows from regularity of the cardinal $\mathfrak{c}$ that $(2^{<\mathfrak{c}})^{|\alpha|}=2^{<\mathfrak{c}}$ for every $\alpha < \mathfrak{c}$. Let now $M_0 \prec H(\theta)$ be such that $A \subset M_0$ and $|M_0| \leq 2^{<\mathfrak{c}}$. Suppose we have constructed $\{M_\alpha: \alpha < \beta \}$ such that for every $\alpha < \beta$ we have $M_\alpha \prec H(\theta)$, $|M_\alpha| \leq 2^{<\mathfrak{c}}$. Then let $M_\alpha \prec H(\theta)$ be such that $M_\beta \cup [M_\beta]^{|\alpha|} \subset M_\alpha$ for every $\beta < \alpha$ and $|M_\alpha| \leq 2^{<\mathfrak{c}}$. Then $\{M_\alpha: \alpha < \mathfrak{c}\}$ is a chain under containment of elementary submodels of $H(\theta)$ and hence it is also an elementary chain, from which it follows that $M=\bigcup_{\alpha < \mathfrak{c}} M_\alpha$ is an elementary submodel of $H(\theta)$.

To see that $M$ is $<\mathfrak{c}$-closed let $\lambda < \mathfrak{c}$ and $\{x_\alpha: \alpha < \lambda \} \subset M$. Then, by regularity of $\mathfrak{c}$ there is $\tau < \mathfrak{c}$ such that $\{x_\alpha: \alpha < \lambda \} \subset M_\tau$. We can certainly assume $\tau > \lambda$. But $[M_\tau]^{|\lambda|} \subset M_{\tau+1}$ and therefore $\{x_\alpha: \alpha < \lambda\} \in M_{\tau+1} \subset M$.
\end{proof}

\begin{lemma}
(Shapirovskii, see \cite{J}, 2.13) Let $\mathcal{U}$ be an open cover for some space $X$. Then there is a discrete $D \subset X$ and a subcover $\mathcal{W} \subset \mathcal{U}$ such that $|\mathcal{W}|=|D|$ and $X= \overline{D} \cup \bigcup \mathcal{W}$.
\end{lemma}

\begin{theorem}
($2^{<\mathfrak{c}}=\mathfrak{c}$) Let $X$ be a space such that $\hat{s}(X) \cdot g(X) \leq \mathfrak{c}$. Then $|X| \leq \mathfrak{c}$.
\end{theorem}

\begin{proof}

Let $M$ be an elementary submodel of a large enough fraction of the universe such that $\{X, \tau\} \subset M$, $\mathfrak{c} \cup \{\mathfrak{c}\} \subset M$, $|M| \leq \mathfrak{c}$ and $M$ is $\lambda$-closed for every $\lambda < \mathfrak{c}$. 

We claim that $X \subset M$. Suppose not and fix $p \in X \setminus M$. We claim that for every $x \in X \cap M$ we can choose an open $U \in M$ such that $x \in U$ and $p \notin U$. Indeed, fix $x \in X \cap M$ and let $\mathcal{V} \in M$ be the set of all open sets $V \subset X$ such that $x \notin \overline{V}$. Then $\mathcal{V}$ covers $X \setminus \{x\}$, so by Shapirovskii's Lemma we can find a discrete $D \in M$ and a subfamily $\mathcal{W} \subset \mathcal{V}$ such that $\mathcal{W} \in M$, $|\mathcal{W}|=|D| \leq \mathfrak{c}$ and $X \setminus \{x\} \subset \overline{D} \cup \bigcup \mathcal{W}$. Now $\mathcal{W} \in M$ and $|\mathcal{W}| \leq \mathfrak{c}$ imply that $\mathcal{W} \subset M$. Notice that, since $D \in M$, also $\overline{D} \in M$ which implies $\overline{D} \subset M$, since $|\overline{D}| \leq \mathfrak{c}$. So $p \notin \overline{D}$ and hence there is $W \in \mathcal{W}$ such that $p \in W$. Let $U=X \setminus \overline{W}$. Then $U \in M$ is a neighbourhood of $x$ such that $p \notin U$.

So for every $x \in X \cap M$ choose $U_x \in M$ such that $p \notin U$. The family $\mathcal{U}=\{U_x : x \in X \cap M \}$ covers $X \cap M$, so, by Shapirovskii's Lemma there is a discrete set $D \subset X \cap M$ and a set $\mathcal{W} \subset \mathcal{U}$ such that $|\mathcal{W}|=|D|<\mathfrak{c}$ with $X \cap M  \subset \overline{D} \cup \bigcup \mathcal{W}$. Since $M$ is $<\mathfrak{c}$-closed we have that $D \in M$ and $\mathcal{W} \in M$, and hence $M \models X \subset \overline{D} \cup \bigcup \mathcal{W}$. Now $p \notin W$ for any $W \in \mathcal{W}$ and $p \notin \overline{D}$, since $\overline{D} \subset X \cap M$, by the same reason as before. But that's a contradiction.
\end{proof}

Can we switch discrete sets with free sequences in the previous theorem? Clearly not, and the one-point compactification of a discrete set is a counterexample. However there are some cases where we can. Let's start by proving a kind of free-sequence version of Shapirovskii's Lemma.

\begin{lemma} \label{shapfree}
Let $X$ be a space such that the closure of every free sequence is Lindel\"of and $\mathcal{U}$ be an open cover for $X$. Then there is a free sequence $F \subset X$ and a subcollection $\mathcal{V} \subset \mathcal{U}$ such that $|\mathcal{V}|=|F|$ and $X = \overline{F} \cup \bigcup \mathcal{V}$.
\end{lemma}

\begin{proof}
Suppose you have constructed, for some ordinal $\beta$, a free sequence $\{x_\alpha: \alpha < \beta \}$ and countable subcollections $\{\mathcal{U}_\alpha: \alpha < \beta \}$ such that $\overline{\{x_\alpha : \alpha < \gamma\}} \subset \bigcup_{\alpha \leq \gamma} \bigcup \mathcal{U}_\alpha$ for every $\gamma < \beta$.

Let $\mathcal{U}_\beta$ be a countable subcollection of $\mathcal{U}$ covering the Lindel\"of subspace $\overline{\{x_\alpha : \alpha < \beta \}}$ and pick a point $x_\beta \in X \setminus \bigcup_{\alpha \leq \beta} \bigcup \mathcal{U}_\beta$. Let $\kappa$ be the least ordinal such that $$\overline{\{x_\alpha: \alpha < \kappa\}} \cup \bigcup_{\alpha < \kappa} \bigcup \mathcal{U}_\alpha=X.$$ Then $\{x_\alpha : \alpha < \kappa\}$ is a free sequence and for $\mathcal{V}=\bigcup_{\alpha < \kappa} \mathcal{U}_\alpha$ we have $|\mathcal{V}|=\kappa$.
\end{proof}

\begin{theorem} \label{thmpsi}
($2^{<\mathfrak{c}}=\mathfrak{c}$) Let $X$ be a Lindel\"of space such that $\psi(X) \leq \mathfrak{c}$ and $\hat{F}(X) \cdot b(X) \leq \mathfrak{c}$. Then $|X| \leq \mathfrak{c}$.
\end{theorem}

\begin{proof}
Let $M$ be a $<\mathfrak{c}$-closed elementary submodel such that $\mathfrak{c} \cup \{\mathfrak{c}\} \subset M$ and $\{X, \tau \} \subset M$. 

\noindent \textbf{Claim:} The closure of every free sequence in $X \cap M$ is Lindel\"of.

\begin{proof}[Proof of Claim]
Let $F \subset X \cap M$ be a free sequence in $X \cap M$ well-ordered in type $\kappa$ (where $\kappa \leq \mathfrak{c}$ because $|M| \leq \mathfrak{c}$). We claim that $F$ is also a free sequence in $X$. Denote by $F_\beta$ the initial segment of $F$ determined by its $\beta$th element. Let $\alpha=\sup \{\beta < \alpha: F_\beta$ is a free sequence in $X$ by the same well-ordering of $F \}$. Then $F_\alpha$ is a free sequence in $X$. If not, there would be some $\beta < \alpha$ such that $x \in \overline{F_\beta} \cap \overline{F_\alpha \setminus F_\beta}$ and $x \notin M$. But $F_\beta$ is a free sequence in $X$ and therefore $|F_\beta| < \mathfrak{c}$. Thus $F_\beta \in M$, and hence $\overline{F_\beta} \in M$, which along with $|\overline{F_\beta}| \leq \mathfrak{c}$ implies that $\overline{F_\beta} \subset M$. So $x \in M$, which is a contradiction. But now $F_{\alpha+1}$ is also a free sequence in $X$, because you can't spoil freeness by adding a single isolated point. Therefore $\alpha=\kappa$, which proves that $F$ is a free sequence in $X$. Proceeding as before we get that $\overline{F} \subset X \cap M$, which proves our claim, since closed subspaces of Lindel\"of spaces are Lindel\"of.
\renewcommand{\qedsymbol}
{$\triangle$}\end{proof}

We claim that $X \subset M$. Suppose not, and let $p \in X \setminus M$. For every $x \in X \cap M$ use $\psi(X) \leq \mathfrak{c}$ to pick a neighbourhood $U_x \in M$ of $x$ such that $p \notin U_x$. Let $\mathcal{U}=\{U_x: x \in X \cap M \}$. By Lemma $\ref{shapfree}$, there are a free sequence $F \subset X \cap M$ and a subcollection $\mathcal{V} \subset \mathcal{U}$ such that $|F|=|\mathcal{V}| < \mathfrak{c}$ with $X \cap M \subset \overline{F} \cup \bigcup \mathcal{V}$. Now $|F|< \mathfrak{c}$, so $F \in M$ and hence $\overline{F} \in M$, which, along with $|\overline{F}| \leq \mathfrak{c}$ implies that $\overline{F} \subset M$. Also, $\mathcal{V} \subset M$ and $|\mathcal{V}| < \mathfrak{c}$ imply that $\mathcal{V} \in M$. Therefore $M \models X \subset \overline{F} \cup \bigcup \mathcal{V}$ and hence there is $V \in \mathcal{V}$ such that $p \in V$, which is a contradiction.
\end{proof}

Pseudocharacter $\leq \kappa$ is not discretely reflexive, unless the space is compact (see \cite{ATW}). The following lemma shows that the pseudocharacter of a space never exceeds its depth.

\begin{lemma} \label{lemchar}
Let $\kappa$ be an infinite cardinal and $X$ be a space where $|\overline{D}| \leq \kappa$ for every discrete $D \subset X$. Then $\psi(X) \leq \kappa$. If in addition $X$ is regular then $\psi(F,X) \leq \kappa$, for every closed $F \subset X$ such that $|F| \leq \kappa$.
\end{lemma}
 
 \begin{proof}
Let $F \subset X$ be a $\kappa$-sized closed set (or a point, if $X$ is not regular). Now let $\mathcal{V}=\{V \subset X: V$ is open and $\overline{V} \cap F=\emptyset \}$. Then $\mathcal{V}$ covers $X \setminus F$ and hence we can find a discrete $D \subset X \setminus F$ and a subcollection $\mathcal{U} \subset \mathcal{V}$ with $|\mathcal{U}|=|D|$ such that $X \setminus F \subset \bigcup \mathcal{U} \cup \overline{D}$. So $(\bigcap_{x \in \overline{D} \setminus F} X \setminus \{x\}) \cap (\bigcap_{U \in \mathcal{U}} X \setminus \overline{U})=F$, which implies that $\psi(F,X) \leq \kappa$.
 \end{proof}
 
 The following corollary is another improvement of Alas' Theorem.
 
 \begin{corollary}
($2^{<\mathfrak{c}}=\mathfrak{c}$) Let $X$ be a Lindel\"of space such that $\hat{F}(X) \cdot g(X) \leq \mathfrak{c}$. Then $|X| \leq \mathfrak {c}$.
 \end{corollary}
 
 \begin{proof}
 This follows from Lemma $\ref{lemchar}$ and Theorem $\ref{thmpsi}$.
 \end{proof}
 
 In the above corollary Lindel\"ofness can be removed, if one assumes the space to be regular.
 
 \begin{theorem} \label{thmreg}
 ($2^{<\mathfrak{c}}=\mathfrak{c}$) Let $X$ be a regular space such that $\hat{F}(X) \leq \mathfrak{c}$ and $|\overline{D}| \leq \mathfrak{c}$ for every discrete $D \subset X$. Then $|X| \leq \mathfrak{c}$.
 \end{theorem}

 \begin{proof}
 Let $M$ be an elementary submodel as before. By Lemma $\ref{lemchar}$ every $\mathfrak{c}$-sized closed subset of $X$ has pseudocharacter $\leq \mathfrak{c}$.

We claim that $X \subset M$. Suppose not and fix $p \in X \setminus M$ and suppose that for some $\beta < \mathfrak{c}$ we have constructed a free sequence $\{x_\alpha: \alpha < \beta \} \subset M$ and open sets $\{U_\alpha: \alpha < \beta \} \subset M$. We have $p \notin \overline{\{x_\alpha : \alpha < \beta \}}$. Now use the claim to choose a sequence $\mathcal{G} \in M$ of open sets such that $|\mathcal{G}| \leq \mathfrak{c}$ and $\overline{\{x_\alpha: \alpha < \beta\}}=\bigcap \mathcal{G}$. We have $\mathcal{G} \subset M$, so we can choose an open set $U_\beta \in M$ with $p \notin U_\beta$ and $\overline{\{x_\alpha: \alpha < \beta\}} \subset U_\beta$. Now use $< \mathfrak{c}$-closed and elementarity to pick $x_\beta \in (X \setminus \bigcup_{\alpha \leq \beta} U_\alpha) \cap M$. Thus $\{x_\alpha : \alpha \leq \mathfrak{c}\}$ is a $\mathfrak{c}$-sized free sequence in $X$, which is a contradiction.
\end{proof}

In Theorem $\ref{thmreg}$ one can safely work in ZFC if free sequences are assumed to be countable. So we have a common framework for Alas' Theorem and Dow's result about compact spaces of countable tightness mentioned in the introduction. We have only one case left to exhaust all relationships between the four cardinal functions we have defined and cardinality.

\begin{theorem}
($2^{<\mathfrak{c}}=\mathfrak{c}$) Let $X$ be a regular space such that $\hat{s}(X) \cdot b(X) \leq \mathfrak{c}$. Then $|X| \leq \mathfrak{c}$.
\end{theorem}

\begin{proof}
Let $F \subset X$. We claim that $\psi(F,X) \leq \mathfrak{c}$. Indeed, for every $x \notin F$ use regularity to choose an open neighbourhood $V_x$ of $x$ such that $\overline{V_x} \cap F=\emptyset$. Then $\{V_x: x \notin F\}$ covers $X \setminus F$, so we can choose a discrete $D \subset X \setminus F$ such that $X \setminus F \subset \bigcup \{\overline{V}_x: x \in D \} \cup \overline{D}$. Now we claim that for every $p \in \overline{D} \setminus F$ we can choose an $E \subset D$ such that $p \in \overline{E}$ and $\overline{E} \cap F=\emptyset$. Indeed, simply use regularity to find an open neighbourhood $U$ of $p$ such that $\overline{U} \cap F=\emptyset$ and set $E=U \cap D$. So $F=\bigcap \{X \setminus \overline{E}: E \subset D$ and $\overline{E} \cap F=\emptyset \} \cap \bigcap \{V_x: x \in D\}$. This implies that $\psi(F,X) \leq \mathfrak{c}$ since $|D| < \mathfrak{c}$ and hence $2^{|D|} \leq \mathfrak{c}$, by the set-theoretic assumption. Now, an argument similar to the proof of Theorem $\ref{thmreg}$ will finish the proof.
\end{proof}

Regularity can be replaced by Lindel\"ofness. We leave the details to the reader. 

\begin{question}
Is there in ZFC a Hausdorff non-regular space such that free sequences are countable (discrete sets are countable), $|\overline{D}| \leq \mathfrak{c}$ for every discrete $D \subset X$ (for every free sequence $F \subset X$) and yet $|X| > \mathfrak{c}$?
\end{question}

\begin{question}
Is there, in some model of set theory, some (compact) regular space $X$ such that every discrete set has size $< \mathfrak{c}$, the closure of every discrete set has size $\leq \mathfrak{c}$ and yet the space has size $>\mathfrak{c}$.
\end{question}

To find a Hausdorff counterexample to the above question, take a model of $\omega_1 < \mathfrak{c} < 2^{\omega_1}$ and let $X=2^{\omega_1}$. Let $\tau=\{U \setminus C: U$ is open in the usual topology on $2^{\omega_1}$ and $|C| \leq \omega_1 \}$. Then every discrete set in $(X, \tau)$ is closed and has size $\omega_1 < \mathfrak{c}$.


\begin{thebibliography}{30}
\bibitem{A} O. Alas, \emph{On closures of discrete subsets}, Q and A in General Topology, Vol. 20 (2002), 85--89.

\bibitem{ATW} O. Alas, V. Tkachuk, R. Wilson, \emph{Closures of discrete sets often reflect global properties}, Topology Proc. 25 Spring (2000).

\bibitem{Ar} A. V. Arhangel'skii, \emph{Structure and classification of topological spaces and cardinal invariants}, Russ. Math. Surv. 33 (1978), 33--96.

\bibitem{D1} A. Dow, \emph{An introduction to applications of elementary submodels to topology}, Topology Proc. (1988), 13(1), 17--72.

\bibitem{D2}A. Dow, \emph{Closures of discrete sets in compact spaces}, Studia Math. Sci. Hung. 42 (2005), 227-234.

\bibitem{E} R. Engelking, \emph{General Topology}, second ed., Sigma Series in Pure Mathematics, no.6, Heldermann Verlag, Berlin, 1989.

\bibitem{Ge1} J. Gerlits, \emph{On a problem of S. Mr\'owka}, Period. Math. Hungar. 4 (1973), 71--80.

\bibitem{Ge2} J. Gerlits, \emph{On a generalization of dyadicity}, Studia Sci. Math. Hungar. 13 (1978), 1--17.

\bibitem{GJS} J. Gerlits, I. Juh\'asz and Z. Szentmikl\'ossy, \emph{Two improvements on Tkachenko's addition theorem}, Comment. Math. Univ. Carolinae, 46, 4 (2005), 705-710.

\bibitem{G1} G. Gruenhage, \emph{Generalized metric spaces} in Handbook of Set Theoretic Topology, edited by K. Kunen and J.E. Vaughan, North-Holland, Amsterdam 1984.

\bibitem{G2} G. Gruenhage, \emph{A note on Gul'ko compact spaces}, Proc. Amer. Math. Soc. 100 (1987), 371-376.

\bibitem{G3} G. Gruenhage, \emph{Covering compacta by discrete and other separated sets}, preprint.

\bibitem{J} I. Juh\'asz, \emph{Cardinal Function in Topology - Ten Years Later}, Mathematical Centre Tracts, 123, Mathematisch Centrum, Amsterdam, 1980.

\bibitem{JM} I. Juh\'asz, J. van Mill, \emph{Covering compacta by discrete subspaces}, Topology Appl. 154 (2007), 283--286.

\bibitem{JS} I. Juh\'asz, Z. Szentmikl\'ossy, \emph{A strengthening of the \v Cech-Pospi\v sil theorem}, preprint.

\bibitem{LS} A. Leiderman, G. Sokolov, \emph{Adequate families of sets and Corson compacts},  Comment. Math. Univ. Carolin.  25 (1984),  no. 2, 233--246.  

\bibitem{R} G.J. Ridderbos, \emph{Power homogeneity in Topology}, Doctoral Thesis, Vrije Universiteit, Amsterdam (2007).

\bibitem{S} S. Spadaro, \emph{Covering by discrete and closed discrete sets}, Topology and its Applications 156 (2009), 721-727.

\bibitem{W} S. Watson, \emph{Locally compact normal spaces in the constructible universe}, Can. J. Math., vol. XXXIV, n.5, 1982, 1091--1096.

\bibitem{Y} N. Yakovlev, \emph{On bicompacta in $\Sigma$-products and related spaces.}, Comment. Math. Univ. Carolin. 21 (1980), no. 2, 263--283. 
\end{thebibliography}
\end{document}